\documentclass{owrart}

\begin{document}

\noindent Oberwolfach Reports {\bf 7} (2010), 2136--2139

\vskip3em

\begin{talk}{Mikhail Belolipetsky}
{Geodesics, volumes and Lehmer's conjecture}
{Belolipetsky, Mikhail}

\vskip2em

\noindent
We first recall a well known relation between Lehmer's conjecture (also known as ``Lehmer's question'') and the Short Geodesic conjecture for arithmetic $2$- and $3$-orbifolds. Here we follow the exposition in~\cite[Chapter~12]{mb_maclachlan_reid}.

Let $P(x)$ be an irreducible monic polynomial with integer coefficients of degree $n$, and let $\theta_1$, \dots, $\theta_n$ denote its roots. The \emph{Mahler measure} of $P$ is defined by
$$ M(P) = \prod_{i=1}^n \mathrm{max}(1, |\theta_i|). $$
If $P(x)$ is a cyclotomic polynomial then its Mahler measure is equal to $1$. Now Lehmer's conjecture says that the measures of all other $P(x)$ are separated from $1$ by an absolute positive constant which is called Lehmer's number:

\medskip

\noindent\textbf{Lehmer's Conjecture.} \emph{There exists $m > 1$ such that $M(P) \ge m$ for all non-cyclotomic $P$.}

\medskip

On the other hand we have the following geometric conjecture.

\medskip

\noindent\textbf{Short Geodesic Conjecture.} \emph{There is a universal positive lower bound for the length of geodesics of arithmetic hyperbolic $2$- and $3$-orbifolds.}

\medskip

To draw the relation between the two conjectures we can argue as follows. Consider a hyperbolic element $\gamma\in\mathrm{PSL}(2,\mathbb{R})$ or $\mathrm{PSL}(2,\mathbb{C})$. We can define its trace and $\gamma$ being \emph{hyperbolic} means that $\mathrm{tr}(\gamma) = u+u^{-1}$ with $|u|>1$ (in case of $\mathrm{PSL}(2,\mathbb{C})$ such elements are often called loxodromic but we will not use this terminology). If $P$ is the minimal polynomial of $u$, then the displacement of $\gamma$ is given by
$$\ell_0(\gamma) = 2\log M(P) \quad \textrm{or} \quad \log M(P)$$
for $2$ and $3$ dimensional cases, respectively. It can be shown that if $\gamma$ is an element of an arithmetic subgroup then $u$ is an algebraic integer and, moreover, its minimal polynomial is not cyclotomic. This leads to a relation between the two conjectures:

\medskip

\noindent\textbf{Corollary 1.} (1) \emph{Lehmer's conjecture implies Short Geodesic conjecture.}\\
(2) \emph{Short Geodesic conjecture implies a special case of Lehmer's conjecture, namely, Lehmer's conjecture for Salem numbers.}

\medskip

Our next goal is to consider the quantitative side of this relation.

The \emph{systole} of a Riemannian manifold $M$, denoted by $\mathrm{Syst}_1(M)$, is the length of a closed geodesic of the shortest length in $M$. This notion can be also generalised to Riemannian orbifolds. It is clear from the definition that systole is a geometric invariant of a manifold (or orbifold). By Borel's theorem there exist only finitely many arithmetic hyperbolic $2$- or $3$-orbifolds of bounded volume~\cite{mb_borel}. Therefore, if the Short Geodesic Conjecture is false then there should exist a sequence of orbifolds $M_i$ such that $$\mathrm{Syst}_1(M_i) \to 0 \quad \textrm{and} \quad \mathrm{Vol}(M_i) \to \infty.$$
Using our current knowledge of Mahler measures of polynomials and some recent results about volumes of arithmetically defined orbifolds we can say more here.

Let $\gamma$ be an element of the group uniformising $M_i$ with the smallest positive displacement and $P_i$ is the minimal polynomial associated to $\gamma$ as above so that $\mathrm{Syst}_1(M_i) = 2\log M(P_i)$ or $\log M(P_i)$ depending on the dimension.

First, let us recall a well known Dobrowolski's bound for the Mahler measure:
\begin{equation}\label{eq1}
\log M(P_i) \ge C_1\left(\frac{\log\log d_i}{\log d_i}\right)^3,
\end{equation}
where $d_i$ is the degree of the polynomial $P_i$ and $C_1>0$ is an explicit constant. To my knowledge, Dobrowolski's result is still the best of its kind except possibly for the value of constant $C_1$. We refer to an excellent survey article by C.~Smyth~\cite{mb_smyth} for more about this topic.

Now we claim that the field of definition $k_i$ of the arithmetic orbifold $M_i$ satisfies
\begin{equation}
\mathrm{deg}(k_i) \ge \frac12 d_i.
\end{equation}
This indeed follows from the results discussed in~\cite[Chapter 12]{mb_maclachlan_reid}.

Finally, we can bring facts together using an important inequality which relates the volume of an arithmetic hyperbolic orbifold with its field of definition:
\begin{equation}\label{eq2}
\mathrm{deg}(k_i) \le C_2\log\mathrm{Vol}(M_i) + C_3.
\end{equation}
This result was first proved by Chinburg and Friedman, in a form stated here it can be found in~\cite{mb_bgls}. We note that the constants in inequalities \eqref{eq1}--\eqref{eq2} can be computed explicitly.

For sufficiently large $x$, $\frac{\log x}{x}$ is a monotonically decreasing function hence for sufficiently small $\mathrm{Syst}_1(M_i)$ we get
\begin{equation}\label{eq3}
\mathrm{Syst}_1(M_i) \ge C_1\left(\frac{\log\log\log \mathrm{Vol}(M_i)^C}{\log\log \mathrm{Vol}(M_i)^C}\right)^3.
\end{equation}

This is a \emph{very slowly decreasing function}! It shows that if the Short Geodesic conjecture is false, it has to be violated by a sequence of orbifolds with extremely fast growing volumes. We also get the following immediate corollary.

\medskip

\noindent\textbf{Corollary 2.} \emph{If $\mathrm{Syst}_1(M_i)\to 0$, then the arithmetic orbifolds $M_i$ are non-commen\-su\-rable and defined over fields of degree going to $\infty$.}

\medskip

Let us now consider the general, not necessarily arithmetic, hyperbolic $n$-manifolds and orbifolds. For dimensions $n\ge4$, in a straight contrast with $n=2$ and $3$, we then have an analogue of Borel's theorem cited above. This result, proved by H.~C.~Wang in~\cite{mb_wang}, implies that for $n\ge 4$ there exist only finitely many hyperbolic $n$-orbifolds of bounded volume. Hence again if we would like to have a sequence of higher dimensional hyperbolic manifolds or orbifolds $M_i$ with $\mathrm{Syst}_1(M_i) \to 0$, we would necessarily have $\mathrm{Vol}(M_i) \to \infty$. In line with the previous discussion and known rigidity properties of higher dimensional hyperbolic manifolds we come to the questions of whether such sequences do exist at all and if yes, then what can we say about the isosystolic properties of corresponding manifolds. The first question was answered for $n=4$ in a remarkable short paper by I.~Agol~\cite{mb_agol}. In a joint work with S.~A.~Thomson we modify and extend Agol's argument which allows us to construct hyperbolic $n$-manifolds with short systole uniformly for all dimensions $n$ and also to answer the second question about their isosystolic properties~\cite{mb_belolipetsky_thomson}. Our main results can be summarised as follows.

\medskip

\noindent\textbf{Theorem 1.}\emph{
\begin{itemize}
\item[(A)] For every $n\geq 2$ and any $\epsilon>0$, there exist compact $n$-dimensional hyperbolic manifolds $M$ with $\mathrm{Syst}_1(M)<\epsilon$.
\item[(B)] For every $n\geq 3$ there exists a positive constant $C_n$ (which depends only on $n$), such that the systole length and volume of the manifolds obtained in the proof of part~(A) satisfy
	$\mathrm{Vol}(M) \geq C_n/\mathrm{Syst}_1(M)^{n-2}$.
\end{itemize}
}

\medskip

Concerning part~(B), we can show that it is possible to achieve that $\mathrm{Vol}(M)$ grows exactly like a polynomial in $1/\mathrm{Syst}_1(M)$, therefore, Theorem~1(B) captures the growth rate of the volume in our construction. This growth can be compared with inequality \eqref{eq3} which says that if a similar phenomenon can occur in arithmetic setting, the volume would have to grow much faster. It is unknown if for $n\geq 4$ there exist hyperbolic $n$-manifolds $M$ with $\mathrm{Syst}_1(M)\to 0$ and $\mathrm{Vol}(M)$ growing slower than a polynomial in $1/\mathrm{Syst}_1(M)$. Let us also remark that an alternative proof of part~(A) of Theorem~1 can be given using the original Agol's construction combined with a recent work of Bergeron, Haglund and Wise~\cite{mb_bhw}.

To conclude the comparison with the low dimensional case and Corollary~2 stated above, here we have

\medskip

\noindent\textbf{Corollary 3.} \emph{When $\mathrm{Syst}_1(M) < \epsilon$, the manifolds $M$ from Theorem~1 are non-arithmetic and all but finitely many of them are non-commensurable to each other.}

\medskip

For the proof of this corollary, its relation to Lehmer's conjecture and some other results we refer to Section~5 of~\cite{mb_belolipetsky_thomson}.

\vskip3em

{\small
\noindent 
Department of Math. Sciences, Durham University, South Rd, Durham, DH1 3LE, UK\\
{\it and}
Sobolev Institute of Mathematics, Koptyuga 4, 630090 Novosibirsk, RUSSIA\\
{\it E-mail address:} mbel@math.nsc.ru
}
\end{talk}

\end{document}